\documentclass{article}
\usepackage{caption}
\usepackage{graphicx}
\usepackage{epstopdf}
\usepackage{fancyhdr}
\usepackage{afterpage}
\usepackage[letterpaper,left=3cm,right=3cm,top=2cm,bottom=2cm,marginparwidth=1.75cm]{geometry}
\usepackage{booktabs}
\usepackage{cite}
\usepackage{amsmath}
\usepackage{mathrsfs}
\usepackage{amsthm,amssymb}
\let\nofiles\relax 
\usepackage{indentfirst}
\usepackage{bm}
\usepackage{url}
\usepackage[colorlinks=true,linkcolor=blue,citecolor=blue,urlcolor=blue]{hyperref}
\usepackage{makecell}
\usepackage{multirow}
\usepackage{slashbox}
\usepackage[ruled,linesnumbered]{algorithm2e}
\usepackage{mathrsfs}
\usepackage{cleveref}
\usepackage[english]{babel}
\usepackage{float}
\usepackage{subfigure}
 
\begin{document}
\thispagestyle{empty} \vspace*{-1.0truecm} \noindent
\vskip 10mm

\begin{center}
\large\bf 3D Surface Reconstruction and Volume Approximation via the meshless methods
\end{center}
\begin{center}
\rm\bf T.Li$^{\textbf{1}}$,
M. Lei$^{\textbf{1,*}}$, James Snead$^{\textbf{2}}$,  C. S. Chen$^{\textbf{2}}$
\end{center}

\begin{center}{
\small\it
\textsuperscript{1}\!\! College of mathematics, Taiyuan University of Technology, Taiyuan 030024, PR China\\
\textsuperscript{2}\!School of Mathematics and Natural Sciences, University of Southern Mississippi, Hattiesburg, MS 39406, USA\\
$*$ Corresponding author. Email address: leimin@tyut.edu.cn
}
\end{center}

\noindent{\small {\large\bf Abstract} \ \
In this paper, we propose several mathematical models for 3D surface reconstruction and volume estimation from a set of scattered cloud data. Three meshless methods including the interpolation-based method by RBF, PDE-based approach by Kansa's method and the Method of Fundamental Solutions are employed and compared.
For the optimal recovery of the surfaces, the selection of free parameters in related PDE models are further studied and analyzed.
Besides, several criteria like Hausdorff distance are employed in above methods instead of the classical parameter $\lambda$ determination strategy, which leads to a more reliable reconstruction
performance. 
Finally, the volume estimation of 3D irregular objects is proposed based on the  optimal reconstructed geometric models via proposed meshless methods. Numerous
numerical examples are presented to demonstrate the effectiveness of the proposed surface reconstruction methods and the volume estimation strategy.

}

\section*{1. Introduction}

A common problem arose from the computer graphics, medical imaging and computer aided design is the reconstruction of a 3D surface defined in terms of point cloud data.
Nowadays, the reconstruction approaches attract more and more researchers' attentions.
These methods are generally divided into explicit methods and implicit methods according to the topology connection constructed from discrete points. 
The implicit surface reconstruction is well-known for its simplicity, where the output points are always considered as the zero-level set of a sampled, evolving 4D function. 
Hoppe et al.\cite{1992Hoppe} constructed a continuous function to reconstruct a 3D surface from a set of unorganized points based on the signed distance function and established a framework for the implicit reconstruction method. 
Dinh et al. \cite{2001Dinh} reconstructed the 3D data cloud points based on Principal Component Analysis (PCA) method and anisotropic basis function without any prior knowledge of topological structure 
as well as point normal vectors. 
Subsequently, the poisson surface reconstruction algorithm was proposed by Kazhdan et al. \cite{2006Kazhdan, 2013Kazhdan},
which fits implicit surfaces by solving poisson equations.
The quality of reconstructed surfaces is further ameliorated through the improved poisson algorithm in 2013. Many other implicit techniques can be referred to \cite{2003Alexa, 1991Muraki, 2006Hornung, 2007Alliez, 2009Huang, 2009Oztireli,liu5710859,Palomar2016Surface,ZENG2022102062}.

Recently, meshless methods have played an important role in surface reconstruction due to their simplicity and effectiveness.
Carr et al. \cite{2001Carr} proposed the interpolation-based meshless radial basis function (RBF) algorithm to reconstruct 3D models, where the signed distance function was considered as the weighted sum of a series of RBFs.
In the calculation process, due to the introduction of auxiliary points in the $\varrho$-width narrow band of the surface, two sets of additional point cloud data relative to the surface need to be employed. 
Therefore, the reconstruction process actually takes 3N points instead of the original N points, which leads to a larger system to be solved.
Besides, the reconstruction performance relies heavily on the choosing of $\varrho$, whose determination is a tough job. 
In order to bypass the employment of additional points and reduce the resultant matrix system,
Ohtake et al. proposed another interpolation-based methods employing multilevel compact RBFs as well as the local quadratic term which was used to approximate the surface in a local vicinity of original points.
Another technical pathway within the implicit approach is the PDE-based method.
Zhao et al. \cite{2001Zhao} adopted the idea of variational method to output a smooth surface based on the level set method. 
A potential method of fundamental solutions was proposed by Tankelevich et al \cite{2009Tankelevich} for surface reconstruction through solving harmonic and biharmonic boundary value problems, in which no normal data were required. 
Improved versions of these classical methods can further be found in \cite{1990Kansa, 2010fili,2016Dou, 1998Fairweather, 1998Golberg}.

In this paper, we focus on three implicit reconstruction approaches based on meshless methods
including:
the interpolation-based method by RBF-MQ (MQ),
PDE-based algorithm by Kansa's method with RBF-MQ (Kansa) and the PDE-based method by Method of Fundamental Solutions (MFS).
Compared with the interpolation methods and PDE models proposed in \cite{Tankelevich2006PotentialFB},
it has been found that PDEs containing free parameters are more powerful in 3D surface reconstruction.
\!By choosing an optimal parameter, the spurious surface generated in the reconstruction procedures can be easily eliminated \cite{liu2020implicit,zheng2020improved,chen2021improved}.
However, the selection of the parameter is still an open problem. 
Apart from the experience formula used in \cite{zheng2020improved}, the Hausdorff distance \cite{1998Cignoni} describing the degree of similarity between two sets of points, \textcolor[rgb]{0.00,0.00,0.00}{as well as another two measurements: the Symmetric Chamfer distance (SCd) and the absolute average distance (AAd), are employed in this paper, which guarantees the precision of reconstructed surfaces and volume estimation.}

The overall context of the remaining paper is given as follows. 
Section 2 mainly introduces the reconstruction process through three approaches. 
In section 3, the selection of free parameters in PDE models is studied and the volume estimation of the reconstructed model is proposed. 
Numerical examples are presented to illustrate the feasibility of 
proposed reconstruction method and the accuracy of volume estimation
strategy in section 4. 
In section 5, some conclusions and suggestions for future research are presented.

\section*{2. Methodologies for 3D surface reconstruction}
\section*{2.1 Interpolation-based method by RBF-MQ (MQ) }
\indent
Assume the points cloud
${\mathbf{x}}_{i}
\in \mathbb{R}^{d}$, $
  i=1, \ldots, N$ are from a manifold (surface) $\mathfrak{M}$
 and satisfy
the following equation
in the implicit surface reconstruction
\begin{equation}\label{eq030}
  u(\mathbf{x})=0,
  x \in \mathfrak{M}.
\end{equation}
The approximation surface 
$\mathfrak{M}^{*}$ can be coincided with the zero contour (isosurface) of $u$, which means $\mathfrak{M}^{*}$ in $\mathbb{R}^{d}$ can be seen as the set
\begin{equation}\label{eq070}
  \mathfrak{M}^{*}=\left\{\mathbf{x}^{*} \in \mathbb{R}^{d}: u(\mathbf{x}^{*})=0\right\}.
\end{equation}
In the meshless method, $u(\mathbf{x})$
can be approximated by a linear combination of basis functions like
\begin{equation}\label{eq020}
  u(\mathbf{x})=\sum_{j=1}^{N} \alpha_{j} \phi\left(\left\|\mathbf{x}-
  \bm{\xi}_{j}\right\|\right), \quad \mathbf{x} \in \mathfrak{M},
\end{equation}
where $\alpha_j$ is the unknown weight and $\phi(||\cdot||)$ is radial basic function.
The $\bm{\xi}_j$ are called the center nodes while $\bf{x}$ are the collocation nodes/point cloud data. Many other widely used RBFs in Table 1 can also be employed in Eq.
(\ref{eq020}).

\vspace{0.5cm}
\begin{table}[ht]
\begin{center}
  \begin{tabular}{|c|c|}
\hline $\phi$ & Formula \\
\hline Cubic & $\phi(r)=r^{3}$ \\
\hline Normalized Multiquadrics & $\phi(r)=\sqrt{1+r^{2} c^{2}}, c>0$ \\
\hline Polyharmonic & $\phi(r)=r^{2 n} \log r, n \geq 1$ in 2D \\
Splines & $\phi(r)=r^{2 n-1}, n \geq 1$ in 3D \\
\hline Gaussian & $\phi(r)=e^{-c r^{2}}, c>0$ \\
\hline
\end{tabular}
\end{center}
\vspace{-0.1cm}
\caption{Classical Radial Basis Functions}
\end{table}
In the surface reconstruction via interpolation-based methods, extra points called off-surface nodes $\bar{\mathbf{x}}$ shown in Figure 1 are always needed to be placed at a proper position away from the boundary to obtain an nontrivial interpolant. 
There are many ways to construct the  $\bar{\mathbf{x}}$ like adding
a small perturbation
$\varrho$ along the normal vectors $
\mathbf{n}_i=
(n_{x,i},n_{y,i},n_{z,i})
$ as follows
\begin{equation}\label{eq044}
\bar{\mathbf{x}}=
{\bf{x}_{i}}
\pm\varrho {\bf{{n}}_{i}}.
\end{equation}
The $\varrho $ is critical in the reconstruction process which is always set to be \cite{fasshauer2007meshfree}
\begin{equation}\label{eq031}
  \varrho=0.01\times d_{max},
\end{equation}
where $d_{max}$ is the maximum value of the edge length of the minimum bounding box covering the $\mathfrak{M}$.
Then the following system can be constructed to determine the unknown coefficient

The final system in matrix form can be given as
\begin{equation}\label{eq040}
  \mathbf{A}\bm{\alpha}=\mathbf{b},
\end{equation}
where $\bm{\alpha}=[\alpha_1,...,\alpha_N]^T$ and
$\bm{A}=[A_{ij}]$ with
\begin{equation}\label{eq041}
\mathbf{A}_{ij}=
\left[
\begin{array}{c}
\phi
(||{\mathbf{x}}_i-\bm{\xi}_j
||)
\\
\phi
(||\bar{{\mathbf{x}}}_i-\bm{\xi}_j
||)
\end{array}
\right],
\quad
\mathbf{b}=
\left[
\begin{array}{c}
\mathbf{0} \\
\mathbf{\pm1}
\end{array}
\right].
\end{equation}
The $\bm{\xi}_j$ are generally taken to be the same with all collocation nodes and then the unknown coefficients $\bm{\alpha}$ can be obtained by solving above system.
The next step we need to construct the meshgrid points $\mathbb{X}$ lying uniformly in the minimum bounding box covering $\mathfrak{M}$, after which one needs to pick out all the points $ \mathbb{X}^{*}$ satisfying
\begin{equation}\label{eq071}
  \mathbb{X}^{*}=\left\{\mathbb{X} \in \mathbb{R}^{3}: u(\mathbb{X})=0\right\}.
\end{equation}
Finally, the MATLAB commands \emph{contour} (\emph{isosurface}) and \emph{patch} are used to obtain the $\mathfrak{M}^{*}$.
The ray tracing technique (RTT) or Dual-Contouring algorithm (DCM) is recommended for higher reconstruction quality.
Numerous methods to estimate $\mathbf{n}_i$ can be found in \cite{LIU2016147} when the points cloud does not includes any information of normal vectors.
It is worth pointing out that the procedures following the determination of the unknown $\bm{\alpha}$ are common steps for various implicit surface reconstruction approaches based on meshless methods. Hence we will not elaborate on this part further.

\subsection*{2.2 PDE-based method by Kansa RBF (Kansa) }
\indent
It can be found that the number of nodes used in above approach increase markedly in system Eq(\ref{eq041}).
The issues like computational efficiency as well as data storage limit its application greatly.
Nowadays more and more researchers focus on 3D reconstruction through solving specified PDEs.
One of the main advantages of using PDEs model is that it is suitable for the reconstruction via incomplete data \cite{liu2020implicit}.
In this section, the following PDEs with boundary conditions need to be solved to reconstruct
$\mathfrak{M}$,
\begin{equation}\label{eq101}
\begin{aligned}
  \mathfrak{L} u (\tilde{\mathbf{x}}) & ={f},
  \quad \tilde{\mathbf{x}} \in \tilde{\mathfrak{M}}, \\
  \mathfrak{B}u(\mathbf{x})&
  =g, \quad \mathbf{x}\in \mathfrak{M}.
\end{aligned}
\end{equation}
where $\tilde{\mathfrak{M}} $ is the inside area of $\mathfrak{M}$ and $\tilde{\mathbf{x}}$ are the inside nodes of the problem domain  \cite{zheng2020improved}.
The PDE model like $\mathfrak{L}=\Delta^2,
\mathfrak{B}=[\mathbf{1},\frac{\partial}{\partial n}]^T$ has been considered.
However, skillful technique needs to be employed to eliminate unwanted spurious surfaces.
These skills like adding some extra nodes in proper positions are always tedious and desirable for complicated 3D surfaces.
In this subsection, 
the model 
$\mathfrak{L}=\Delta -\lambda$ with $\mathfrak{B}=\mathbf{1}$ is considered since the unwanted spurious surfaces can be easily eliminated by choosing a proper $\lambda$
as well as it is effective for a wide range of $\lambda$ \cite{liu2020implicit}. More details about choosing proper $\lambda$ will be discussed in following sections.
When we set ${f}=g=1$, the Eq.(\ref{eq070}) are updated as
\begin{equation}\label{eq170}
  \mathfrak{M}^{*}=\left\{\mathbf{x}^{*} \in \mathbb{R}^{3}: u(\mathbf{x}^{*})=1\right\}.
\end{equation}
Then the Kansa's collocation method is employed here due to its high accuracy in solving  high dimensional problems in irregular area.
The Normalized multiquadric  (NMQ) with free shape parameter $c$ is used in this section. 
Based on the interpolation in (\ref{eq020}), we get the matrix system like 
Eq.(\ref{eq041}) corresponding to PDE (\ref{eq101}), where
$\mathbf{A}$ and $\mathbf{b}$ becomes
\begin{equation}\label{eq0401}
\mathbf{A}_{ij}=
\left[
\begin{array}{c}
\mathfrak{L} \phi
(||\tilde{\mathbf{x}}_i-\bm{\xi}_j
||),
\\
\phi
(||{\mathbf{x}}_i-\bm{\xi}_j
||)
\end{array}
\right],
\quad
\mathbf{b}=
\left[
\begin{array}{c}
\mathbf{1} \\
\mathbf{1}
\end{array}
\right].
\end{equation}
By solving above system, we can obtain the unknown coefficients
$\bm{\alpha}$ and the subsequent reconstruction steps are similar with the ones in above section.
Note that when the number of $\tilde{\mathbf{x}}_i$ is less than $2N$, a smaller system can be obtained than the one of interpolation-based methods in Eq.(\ref{eq041}).


\subsection*{2.3 PDE-based method by MFS (MFS) }
The computational cost of above method is still expensive since a set of inner nodes $\tilde{\mathbf{x}}$ is required. In order to further simplify the reconstruction procedures and reduce the computational time,
the PDE model and related numerical methods need to be chosen carefully.
Inspired by the PDE model (Model I) \cite{zheng2020improved} where
\begin{align}\label{eq201}
  & \left
  (\Delta-\lambda^{2}\right)^2 u(\tilde{\mathbf{x}})=0, \quad\tilde{\mathbf{x}} \in \tilde{\mathfrak{M}},\\
  & u(\mathbf{x})=1, \quad\mathbf{x} \in \mathfrak{M},\\
  &\frac{\partial}{\partial \mathbf{n}}u(\mathbf{x})=g(x), \quad\mathbf{x}\in \mathfrak{M},
\end{align}
we consider an improved PDE model (Model II) , where Eq.(\ref{eq201}) becomes
It can be found that the reconstruction performance of Model-II is more stable than Model-I for a wide range of $\lambda$, which will be detailed in examples in Section 4.
The MFS is employed to solve above PDEs due to its simplicity and effectiveness.
As a boundary type meshless method, the discretization only needs to ba carried on the boundaries.
Compared with the traditional global type meshless methods, a smaller linear system will be solved.
The combination of homogeneous PDE model and meshless MFS here will result in a high efficient reconstruction process as the related fundamental solutions are known.
In MFS, the numerical solution of Eq.(17) can be  approximated via a  combination of fundamental solutions as
\begin{equation}\label{eq023}
  u(x, y, z)=\sum_{i=1}^{N}
  \sum_{j=1}^{2} \alpha_{i,j} G_{j}(r, \lambda), \quad(x, y, z) \in \mathfrak{M}.
\end{equation}
where $G_{j}(r, \lambda)$ are the fundamental solutions taking the form 
{
\textcolor[rgb]{1.00,0.00,0.00}{G1 is correct???}
}
\\
\begin{equation}\label{eq0244}
  G_1(r, \lambda)=\frac{e^{-\lambda r}}{4 \pi r}, 
  \quad\quad
G_2(r, \lambda)
=\left\{\begin{array}{l}
\frac{e^{-\lambda r}}{8\pi\lambda} , \quad\text{Model-I},\\
\frac{e^{-\lambda r}-1}{r} , \quad\text{Model-II}.
\end{array}\right.  
\end{equation}
{
Substituting Eq.(\ref{eq0244}) into Eq.(13)-(15) in MFS, we get 
\begin{equation}\label{eq0300}
  \begin{aligned}
&
\sum_{i=1}^{2}A_{1i} \mathbf{a}_{.,i}=\mathbf{1}, \\
&\sum_{i=1}^{2}A_{2i} \mathbf{a}_{.,i}=\mathbf{g},
\end{aligned}
\end{equation}
where $\mathbf{a}_{.,\mathfrak{v}}=[
{a}_{1,\mathfrak{v}},
{a}_{2,\mathfrak{v}},...,
{a}_{N,\mathfrak{v}}]^T,
\mathfrak{v}=1,2$. $\mathbf{1}=[1,1, \cdots, 1]^{T}_{N\times 1}, \mathbf{g}=\left[g\left(x_{1}, y_{1}, z_{1}\right), \cdots, g\left(x_{N}, y_{N}, z_{N}\right)\right]^{T}$ and
$$
\begin{aligned}
A_{1\mathfrak{v}} &=\left[G_{\mathfrak{v}}\left(r_{i j}, \lambda\right)\right]_{1 \leq i, j \leq n}, &
A_{2\mathfrak{v}} &=\left[\frac{\partial}{\partial n} G_{\mathfrak{v}}\left(r_{i j}, \lambda\right)\right]_{1 \leq i, j \leq n}.
\end{aligned}
$$
To accelerate the solving procedure further, $\mathbf{g}$ in Eq.(\ref{eq0300}) can be chosen skillful such that $\mathbf{\alpha}_{.,1}=\mathbf{0}$. Thus Eq.(\ref{eq0300}) becomes
\begin{equation}\label{eq031x}
  A_{12} \mathbf{\alpha}_{.,2}=\mathbf{1},
\end{equation}
by which we can obtain the unknown coefficient $\mathbf{\alpha}_{.,2}$ easily.
Finally the approximation in Eq.(\ref{eq023}) can be simplified further as
\begin{equation}\label{eq0233}
  u(x, y, z)=\sum_{i=1}^{N}
   \alpha_{i,2} G_{2}(r, \lambda),
\end{equation}
which can be used to evaluate the function values on evaluation nodes $\mathbb{X}$ in 
Eq.(\ref{eq071}) to obtain $\mathbb{X}^*$.
The novelty of this method is that normal vectors are no longer needed, which greatly expands its ability in reconstructing complex models,
especially for the ones whose normal vectors are expensive to obtain. 

{\bf{Remark}}:
There is only one free parameter $\lambda$ whose value is critical for the quality of reconstruction.
Based on various numerical examples, for manifold the points cloud data with different $R_{min}$, the product $\lambda R_{min}$ should lie in an interval like \cite{zheng2020improved}
\begin{equation}
0.020\leq\lambda R_{min}\leq0.084. 
\end{equation}
The $R_{min}$ denotes the average smallest distance between the points cloud data $\mathbf{x}_i$ and centers $\bm{\xi}_j$, which can be calculated by 
\begin{equation}\label{eq3001}
  R_{min}=\frac{1}{N}
  \sum_{i=1}^{N}
  min_{1\leq j\leq N}((\mathbf{x}_i,\bm{\xi}_j)).
\end{equation}
The centers and collocation nodes can always taken to be the same since the fundamental solution $G_2$ in Eq.(\ref{eq0233}) is no longer singular.

\section*{3.  Parameter selection and Volume estimation}

\subsection*{3.1 Parameter selection }
Compared with the shape parameter $c$ used in RBF-MQ in Sections 2.1 and 2.2, the $\lambda$ in PDE models in Sections 2.2 and 2.3 plays a more important role in surface reconstruction as well as volume estimation. 
The straightforward strategy way of choosing a suitable $\lambda$ is the experience formula given in 
Eq.(\ref{eq3001}).
However, it will give a range of $\lambda$ and we have no idea which one is the best.
Therefore, we introduce 
the criteria named Hausdorff distance
describing the degree of similarity between two different sets of points.
It gives the maximum value from one set of points to another one, i.e.
\begin{equation}\label{eq321}
  Hd(\tilde{A}, \tilde{B})=\max (h(\tilde{A}, \tilde{B}), h(\tilde{B}, \tilde{A})), 
\end{equation}
where $ h(\tilde{A}, \tilde{B})=\max _{a \in \tilde{A}}\left\{\min _{b \in \tilde{B}}\|a-b\|\right\} $, 
$ h(\tilde{B},\tilde{A})=\max _{b \in \tilde{B}}\left\{\min _{a \in \tilde{A}}\|b-a\|\right\} $.
Here we take    $\tilde{A}=\mathbf{x}_i$ and $ \tilde{B}=\mathbb{X}^{*}$.
\textcolor[rgb]{0.00,0.00,0.00}{
One can also use the absolute average distance (AAd) \cite{cruz2021deepcsr} and Symmetric Chamfer distance (SCd)\cite{park2019deepsdf} defined as}
\textcolor[rgb]{0.00,0.00,0.00}{
\begin{equation} SCd\left(\tilde{A}, \tilde{B}\right)=\frac{1}{m_{\tilde{A}}} \sum_{x \in \tilde{A}} \min _{y \in \tilde{B}}\|x-y\|_2^2+\frac{1}{n_{\tilde{B}} }\sum_{y \in \tilde{B}} \min _{x \in \tilde{A}}\|y-x\|_2^2.
\end{equation}
}
More details of the relationship and difference between above criteria can be found in the following table\\
\begin{table}[h!]
  \centering
  \captionsetup{skip=2pt}
  \caption{The criteria and main difference}
  \begin{tabular}{|c|m{6cm}|m{6cm}|}
    \hline
    \textbf{Method} & \textbf{Description} & \textbf{Key Differentiation} \\
    \hline
    \text{SCd} & Symmetric Chamfer distance calculates the sum of the average squared distances from predicted surface points to the nearest ground truth points, and vice versa. & SCd focuses on the sum of squared distances, highlighting larger distance differences, as squaring amplifies larger distance disparities. Moreover, SCd is more sensitive to outliers, as squared distances amplify larger distance differences. \\
    \hline
    \text{AAd} & AAd measures the mean absolute nearest-neighbor distances between predicted and ground-truth surface points, considering both directions—predicted to ground-truth and ground-truth to predicted. & AAd focuses on the direct mean of distances, emphasizing actual physical distances. Besides, AAd is more intuitive and easier to interpret, as it directly reflects the average distance. \\
    \hline
    \text{Hd-K} & Hausdorff distance identifies the maximum distance from any point on one surface to the closest point on the other surface. & The Hausdorff distance-K is also considered since it can reduce the impact of outliers by considering the K-th percentile of the nearest-neighbor distances in both directions. One can set $K=100$ to use classical Hausdorff distance as this paper does. \cite{gopinath2023learning,walluscheck2023mr}. \\
    \hline
  \end{tabular}
\end{table}

\hspace{2cm}
\\
\textcolor[rgb]{0.00,0.00,0.00}{
Therefore, the optimal $\lambda^{*}$ satisfies
\begin{equation} \lambda^{*}=\underset{\lambda \in \mathbb{R}}{argmin} 
 {\mathfrak{X}d(\mathbf{x}_i, \mathbb{X}^{*})}.
\end{equation}
where $\mathfrak{X}d$ can be $Hd,SCd,AAd$.
}

\subsection*{3.2 {Volume Approximation} }
Based on the optimal reconstructed surface
via above methods with $\lambda^{*}$,
the volume of reconstructed 3D models can be approximated by
\begin{equation}\label{eq080}
  \text { V }=\frac{\text {  $N_{
  {
  \bar{\mathbb{X}^{}}
  }
  }$ }}{\text { 
 $\hat{N}^3$
 }} \cdot \text { $V_{{\mathbb{X}}}$
 ,}
\end{equation}
where $N_{
{
  \bar{\mathbb{X}^{}}
  }
}$ 
\textcolor{black}{
\begin{equation}
\kappa=\left\{\begin{array}{l}
0, \text{MQ} \\
1, \text{Kansa/MFS}.
\end{array}\right.
\end{equation}
}
$\hat{N}^3$ is the total number of 
the uniform nodes $\{{\mathbb{X}}_i\}_{i=1}^{\hat{N}^3}$ in Eq.(\ref{eq071}) in the bounding box
while
$V_{{\mathbb{X}}}$ is the 
volume of bounding box.
The pseudo-code of 3D volume estimation based on interpolation-based method via RBF-MQ is given as an example:
\par
\vspace{3mm}
\begin{algorithm}[H]
\renewcommand{\thealgocf}{}
    \caption{}\label{alg:r2p}
    \KwIn{points cloud data $\mathbf{x}$\;}
       \textbf{Construct} the $\bar{\mathbf{x}}$ by Eq.(\ref{eq044})\;
  \textbf{Solve} the system Eq.(\ref{eq040}) to obtain the unknown coefficients $\bm{\alpha}$\;
  \textbf{Generate} the minimum bounding box containing the $\mathfrak{M}$ and generate $\hat{N}^3$ uniformly distributed nodes 
  ${\mathbb{X}}$ lying in the bounding box\;
  \textbf{Evaluate} the $u({\mathbb{X}})$ and then  use the isosurface command in matlab to obtain the ${\mathfrak{M}}^{*}$ based on Eq.(\ref{eq071}).\;
\textbf{Count} the $N_{
{
  \bar{\mathbb{X}^{}}
  }
}$ and use the formula Eq.(\ref{eq080}) to approximate the volume of the object\;
    \KwOut{the volume of the 3D object\;}
\end{algorithm}

\section*{4. Numerical examples}
In this section, several numerical examples are given to demonstrate the comparison among the detailed reconstruction methods in Section 2 and the 
reliability of proposed volume estimation strategy in Section 3.

\subsection*{4.1 Example 1}
Consider the PDE Model \uppercase\expandafter{\romannumeral1} and \uppercase\expandafter{\romannumeral2} in Section 2.3 to reconstruct the dragon model with 22998 point data respectively. 
To facilitate the comparison of their reconstruction stability with respect to $\lambda$, only the PDE-based method by MFS is employed in this section.
In Figure \ref{Fig.03(a)} ,
the teeth of the dragon are connected when we choose $\lambda=5.83$ based on the experienced formula in Eq.$(23)$\cite{zheng2020improved}. If we increase $\lambda$ to 6.03, the connection part in Figure \ref{Fig.3(b)_1} almost disappears while the region upside the front-legs occurs tiny holes. 
When the $\lambda$ increases a bit to 6.13 in Figure \ref{Fig.3(c)_1} , the spurious connection part disappears totally while the hole becomes larger. 
It is not easy to find a suitable value as
the PDE model I is sensitive to $\lambda$.
We have to pay more time and attention to find the optimal $\lambda$ since the suitable interval is so narrow in some cases.

For the improved Model \uppercase\expandafter{\romannumeral2}, the reconstructed dragon is given in Figure  \ref{Fig.3(d)_1} to \ref{Fig.03(e)}, from which we find the performance is acceptable for a large range of $\lambda\in[21,140]$. 
The spurious part between the teeth vanishes entirely without bringing any holes in the body.
When the $\lambda$ is increased further to $141$ in \ref{Fig.3(f)_1}, the dragon can also be reconstructed well but the surface of dragon's neck becomes non-smooth.
The PDE model II is obviously  more stable and it is then employed in the following examples in MFS without extra mentioning.

\begin{figure}[htbp]
\centering
    \subfigure[(a) Model \uppercase\expandafter{\romannumeral1}: $\lambda=5.83$]{
	\includegraphics[width=6.65cm]{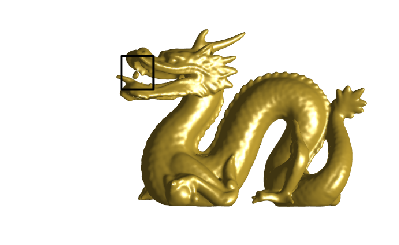} 
    \label{Fig.03(a)}
    }	   
    \hspace{2mm}
 \subfigure[(d) Model \uppercase\expandafter{\romannumeral2}: $\lambda=17$]{
	\includegraphics[width=6.5cm]{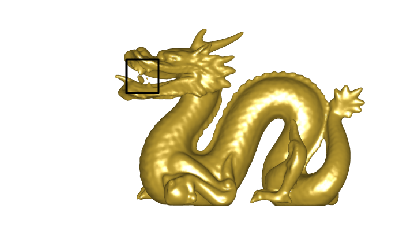} 
    \label{Fig.3(b)_1}
    }    
	\hspace{2mm}
 \subfigure[(b) Model \uppercase\expandafter{\romannumeral1}: $\lambda=6.03$]{
	\includegraphics[width=6.65cm]{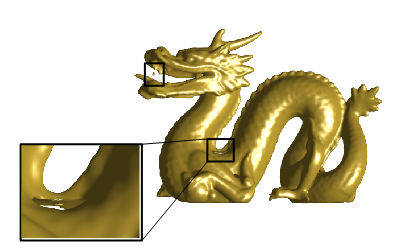} 
    \label{Fig.3(c)_1}
	}
     \hspace{2mm}
 \subfigure[(e) Model \uppercase\expandafter{\romannumeral2}: $\lambda=21,\lambda=45,\lambda=140$]{
	\includegraphics[width=6.5cm]{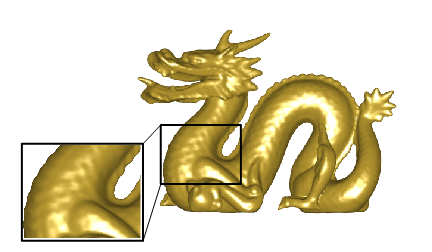} 
    \label{Fig.3(d)_1}
	}
    \hspace{2mm}
 \subfigure[(c) Model \uppercase\expandafter{\romannumeral1}: $\lambda=6.13$]{
	\includegraphics[width=6.55cm]{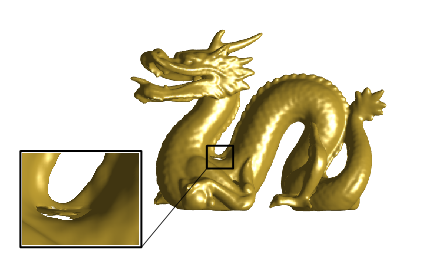} 
    \label{Fig.03(e)}
	}
    \hspace{2mm}
	\subfigure[(f) Model \uppercase\expandafter{\romannumeral2}: $\lambda=141$]{
	\includegraphics[width=6.5cm]{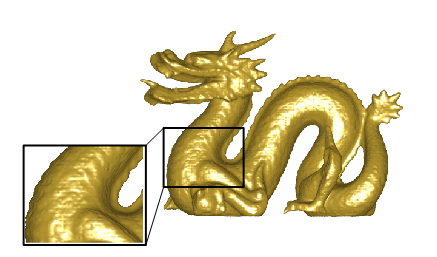} 
    \label{Fig.3(f)_1}
	}
	\caption{The reconstruction comparisons among PDE Model \uppercase\expandafter{\romannumeral1} and \uppercase\expandafter{\romannumeral2} with different $\lambda$ in MFS.}
\label{Figure 03}
\end {figure}


\newpage
\subsection*{4.2 Example 2}
Consider the reconstruction and volume estimation of the 3D bumpy sphere generated by 10000 points through 
\begin{equation}
\mathfrak{M}=\{(x, y, z): x=\rho \sin \phi \cos \theta, y=\rho \sin \phi \sin \theta, z=\rho \cos \phi, \quad 0 \leq \theta \leq 2 \pi, \quad 0 \leq \phi \leq \pi\},
\end{equation}
where $\rho(\phi, \theta)=1+\frac{1}{5} \sin (6 \theta) \sin (6 \phi)$, after which the exact volume can be calculated by
\begin{equation}
    \int_0^{2 \pi} \int_0^\pi\left(\int_0^{1+\sin (6 \theta) \sin (6 \phi) / 5} \rho^2 \sin (\phi) d \rho\right) d \phi d \theta=4.3153.
\end{equation}
The optimal $\lambda^{*}$ in the PDE model is chosen to minimize the \textcolor[rgb]{0.00,0.00,0.00}{ $\mathfrak{X}d$ value.}
The actually optimal value for MFS is $ \lambda^{*} = 9.5729 $ with the minimum $Hd_{min}=0.0735$, 
\textcolor[rgb]{0.00,0.00,0.00}{
where the other optimal $\lambda^{*}$ found by SCd and AAd can be obtained similarly, i.e. $\lambda^{*}=11.8657$,$\lambda^{*}=10.4051$, respectively}.
To facilitate comparison, these actual values of Hd, SCd, and AAd were scaled to the range [0, 1] in Figure \ref{Figure 301} by multiplying an appropriate constant. 
Note that this scaling strategy is also applied to subsequent similar images, unless otherwise stated.
In Figure \ref{Figure 301}, it can be found that although there are slight differences in the optimal $\lambda^{*}$ found via different criteria, their final calculated volume of the object in Tabel \ref{Table 2} remains highly consistent due to the stability of the proposed reconstruction PDE model II in MFS.

More details can be found in Table \ref{Table 2} which also shows the 
volume estimation by Kansa and MQ when $\hat{N}$ = 50, 100,150 respectively.
Compared with the volume $V=4.3164$ given by MeshLab with accuracy $10^{-3}$, it can also be found that the proposed volume estimation performs better in most parameter cases.
Note that for the Kansa method, we choose $ \lambda = 3800$ with $Hd_{min}=0.0748$ when using the radial basis multiquadrics.
The parameter  $ c = 300 $ is used based on experience in radial basis multiquadrics both in the 
interpolation-based method as well as
PDE-based method by Kansa in all the numerical examples.
A detailed view of error distribution is shown in the insets of the Figure \ref{Figure Ex2_1}, where the main differences are emphasised in the red circle.
In Figure \ref{Figure 4_1}, the volumes obtained using different methods and different values of $\hat{N}$ are highly consistent, varying only at the third decimal place.

\textbf{Remark}.
It is worth noting that the $\lambda^{*}$ in MFS is determined by employing a for-loop to choose minimum  metrics such as Hd, SCd, and AAd. The $\lambda^{*}$ significantly impacts the reconstructed surface's quality, such as the watertightness of the generated surface and the presence of any extraneous pseudo-surfaces. 
For simplicity, the optimal $\lambda^{*}$ been found in the case $\hat{N}=50$ in MFS is then directly applied to cases 
$\hat{N}=100,150$ since different 
$\hat{N}$ means varying densities of bounding box grid points, which only leads to a slight variations in  the reconstructed surfaces’s skin textures having a negligible impact on the final volume.

The selection procedure of $\lambda$ in the PDE-based method by Kansa is similar to that in the PDE-based method by MFS and will not be elaborated upon here. Instead, we provide the minimum  Hd related to optimal $\lambda^{*}$ and the corresponding volume. Additionally, the Hd metric and volume for the respective Interpolation-based method by MQ are also presented.
\vspace{3mm}

\begin{table}[htbp]
\vspace{3mm}
  \centering
  \captionsetup{skip=1pt} 
  \caption{The volume of the bumpy model with optimal $\lambda^{*}$ given in Figure \ref{Figure 301} and their related Hd, SCd and AAd obtained by MFS in Example 2}
    \begin{tabular}{c|c|c|c}
    \hline 
    \makebox[0.08\textwidth][c]{$\hat{N}$} &
    \makebox[0.28\textwidth][c]{Volume (Hd, Absolute)} &
    \makebox[0.28\textwidth][c]{Volume (SCd, Absolute)} &
    \makebox[0.28\textwidth][c]{Volume (AAd, Absolute)}\\
    \hline 
    50 & 4.3196 (0.0735, 4.3125e-3) & 4.3207 (8.0118e-04, 5.3950e-03) & 4.3216 (1.8176e-02, 6.2516e-03)\\
   100 &  4.3170 (0.0744, 1.7262e-3) & 4.3174 (5.2178e-04, 2.0874e-03) & 	4.3184	(1.6896e-02, 2.9874e-03)\\
   150 & 4.3172 (0.0735, 1.8935e-3)  &4.3174(4.7171e-04, 2.0474e-03)	& 	4.3182 (1.7440e-02, 2.8574e-03)\\
    \hline  
    \end{tabular}%
  \label{Table 2}%
\end{table}%

\begin{table}[htbp]
\vspace{3mm}
  \centering
  \captionsetup{skip=1pt} 
  \caption{The volumes and related Hd obtained by Kansa and MQ in Example 2}
    \begin{tabular}{c|c|c|c}
    \hline 
    \makebox[0.18\textwidth][c]{Method} &
    \makebox[0.15\textwidth][c]{$\hat{N}$} &
    \makebox[0.28\textwidth][c]{Volume} &
    \makebox[0.28\textwidth][c]{Hd (Absolute Error)} \\
    \hline \multirow{3}[0]{*}{Kansa} & 50 & 4.3164 &  0.0748 ($1.1661e-3$) \\
                                  & 100 & 4.3151 & 0.0758 ($2.0454e-4$) \\
                                  & 150 & 4.3154 & 0.0821 ($1.2432e-4$) \\ 
    \hline \multirow{3}[0]{*}{MQ} & 50 & 4.3161 &  0.0751 ($7.8469e-4$) \\
                                  & 100 & 4.3148  & 0.0785 ($4.6674e-4$) \\
                                  & 150 &4.3154  & 0.0791 ($8.5475e-5$) \\
    \hline  
    \end{tabular}%
  \label{Table 22}%
\end{table}%

\begin{figure}[htbp]
\centering 
\includegraphics[height=7.8cm,width=9cm]{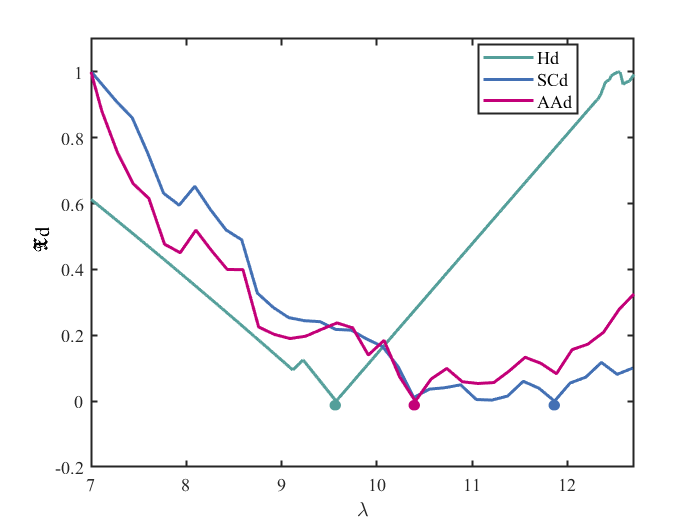}
\captionsetup{skip=1pt} 
\caption{The optimal $\lambda^{*}$ found by different criteria measure when $\hat{N}$ = 50 via PDE-based method by MFS in Example 2.}
\label{Figure 301}
\end{figure}

\newpage
\newgeometry{bottom=2cm}
{
\begin{figure}[htbp]
\centering
    \subfigure[(a) PDE-based method by MFS]{
	\includegraphics[width=12cm]{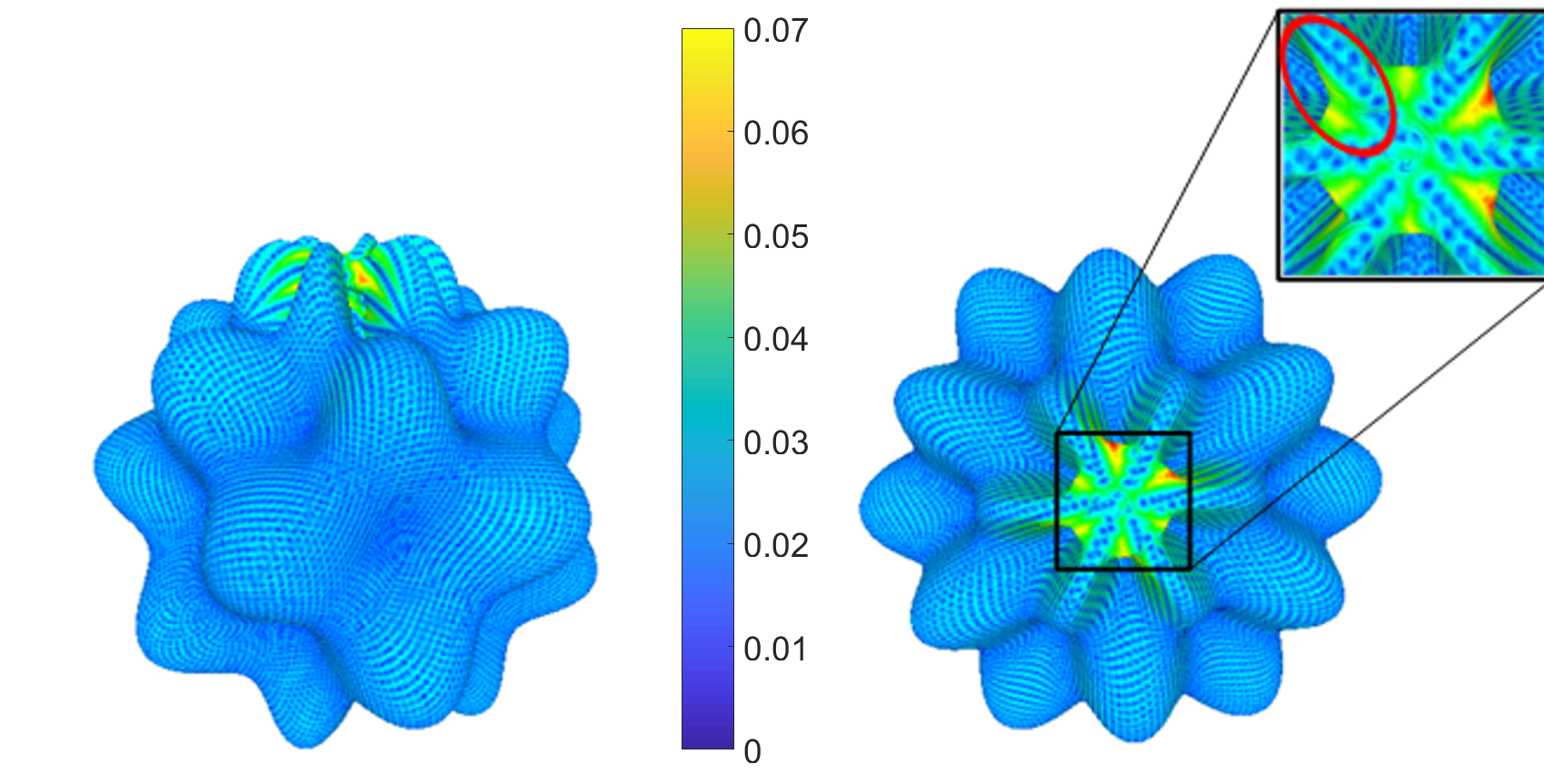} 
    \label{Fig.3(a)}
    }	   
    \hspace{2mm}
	\subfigure[(b) PDE-based method by Kansa]{
	\includegraphics[width=12cm]{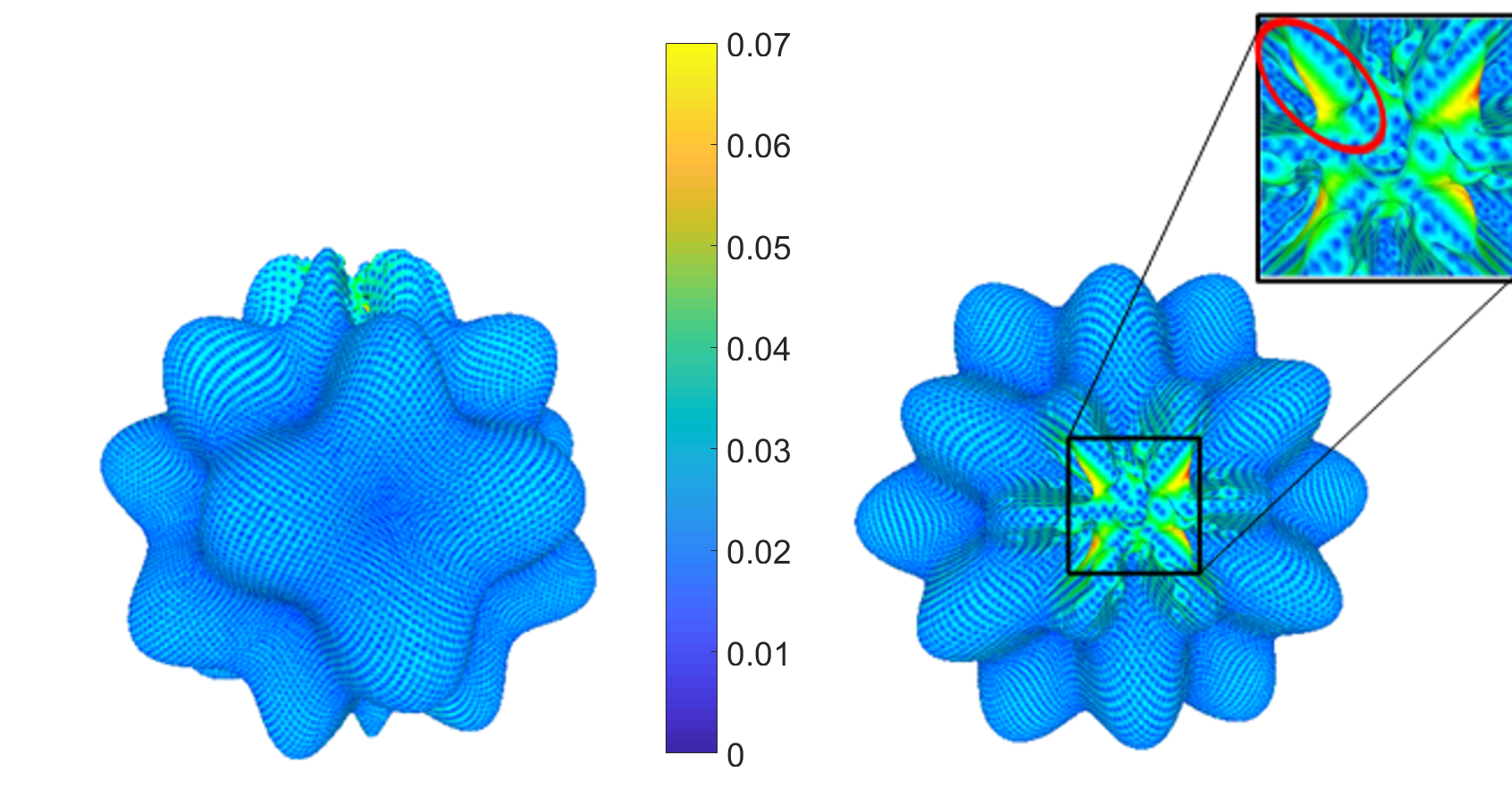} 
    \label{Fig.3(ba)}
    }    
	\hspace{2mm}
	\subfigure[(c) Interpolation-based method by RBF-MQ]{
	\includegraphics[width=11.5cm]{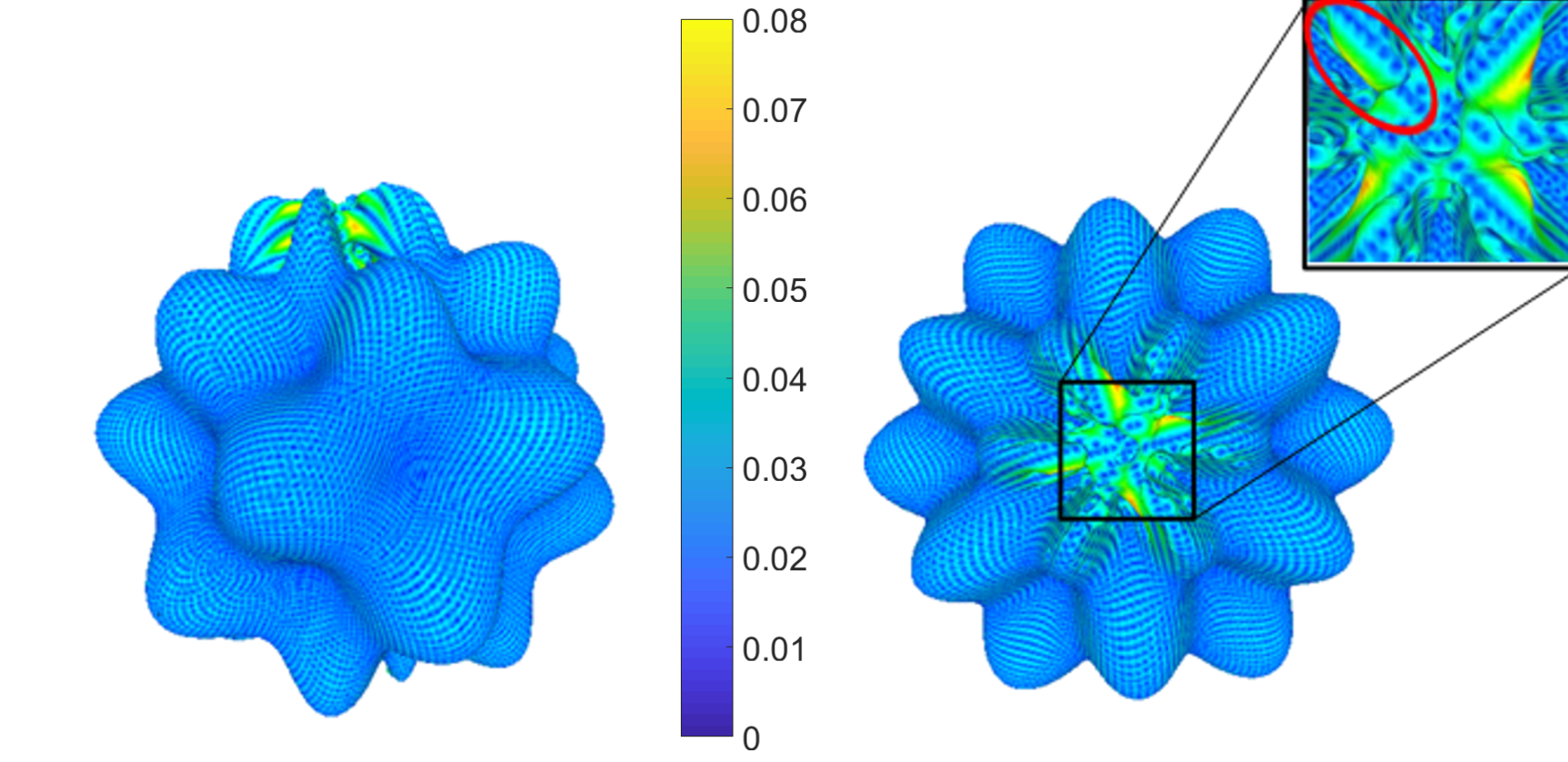} 
    \label{Fig.3(ccc)}
	}
 	\caption{The error distribution of surface reconstruction in Example 2.
  {
  \textcolor[rgb]{1.00,0.00,0.00}{
  (Plotted by MeshLab)
  }
  }
  }
\label{Figure Ex2_1}
\end {figure}
}

\newpage
\vspace{0.7cm}
\begin{figure}[htbp]
\centering 
\includegraphics[height=7.8cm,width=12cm]{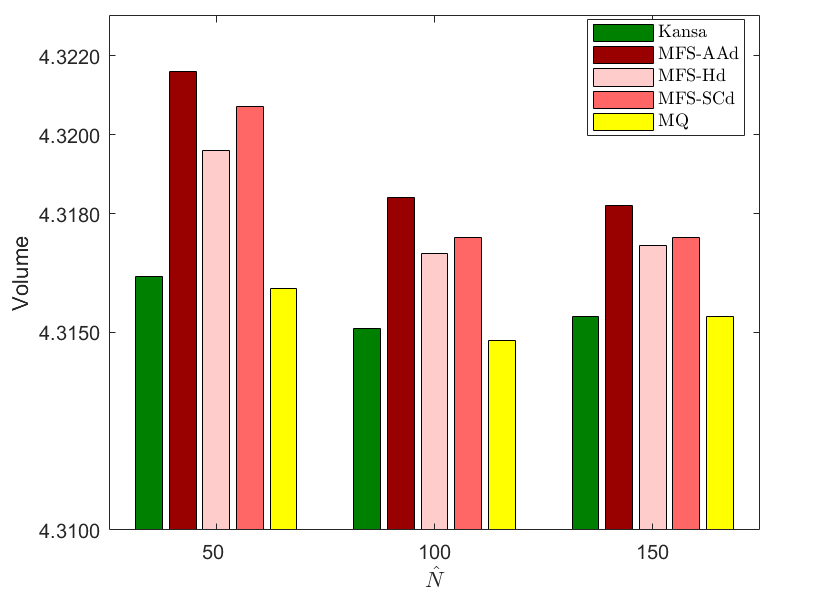}
\caption{Volume under different methods and criteria for Example 2.}
\label{Figure 4_1}
\end{figure}

\subsection*{4.3 Example 3}
In this example, 6967 points data are employed to reconstruct the Stanford Bunny. 
Figure \ref{Figure 5} shows the optimal $\lambda^{*}$ found by Hd,SCd and AAd in the MFS.
For example, the Hd reaches the minimum value when $ \lambda^{*} = 13.7143 $ while
\textcolor[rgb]{0.00,0.00,0.00}{the $\lambda^{*}$  found via SCd and AAD are $\lambda^{*}=6.7931$,$\lambda^{*}=6.1266$, respectively}.
Different criteria may identify inconsistent $\lambda^{*}$, but they yield consistent volumes, as shown in Table 5. It is evident that the magnitude of volume fluctuation is on the order of $10^{-3}$.

For the Kansa method, we choose $ \lambda^{*} = 2052$ which minimizes the Hd as $0.0083$.
The volumes of Bunny by different methods are listed in Table \ref{Table 3} where we can find the volumes are maintained at about 0.754 with little oscillations.
 In Figure \ref{Figure 5_1}, these algorithms produce remarkably consistent volume estimates across different criteria, with variations generally within the range of 
$10^{-3}$. The estimation is almost consistent with the results given in Meshlab with 0.751. While there are slight deviations for specific values of 
$\hat{N}$, the overall trend remains stable, reinforcing the reliability of these computational approaches in volume estimation tasks. 
Figure \ref{Figure Ex3_1} shows the
error distribution of the surface reconstruction by MFS, Kansa and MQ respectively when $\hat{N}$ = 150 in different view. 
The color map in the middle illustrates the related errors measured by Hd in the software Meshlab.
It shows that the main different detailed in the insets comes from the reconstruction of the Bunny's ears.

\vspace{1cm} 
\begin{table}[htbp]
  \centering
  \caption{The volume of the Bunny model with $\lambda^{*}$ given in Figure \ref{Figure 5} and related Hd, SCd and AAd obtained by MFS in Example 3}
    \begin{tabular}{c|c|c|c}
    \hline 
    \makebox[0.1\textwidth][c]{$\hat{N}$} &
    \makebox[0.26\textwidth][c]{Volume (Hd)} &
    \makebox[0.26\textwidth][c]{Volume (SCd)}&
    \makebox[0.26\textwidth][c]{Volume (AAd)} \\
    \hline 
   50    & 0.7539 (0.0076) &0.7547 ($1.9252e-06$) & 0.7549 ($1.0007e-03$)\\
   100   &  0.7546 (0.0084) &0.7557 ($8.4059e-07$) &0.7558 ($5.6838e-04$)
\\
   150   & 0.7546 (0.0087) & 0.7557 ($6.4044e-07$)& 0.7558 	($5.1691e-04$)\\
    \hline  
    \end{tabular}%
  \label{Table 3}%
\end{table}%

\begin{figure}[htbp]
\centering 
\includegraphics[height=7.8cm,width=10cm]{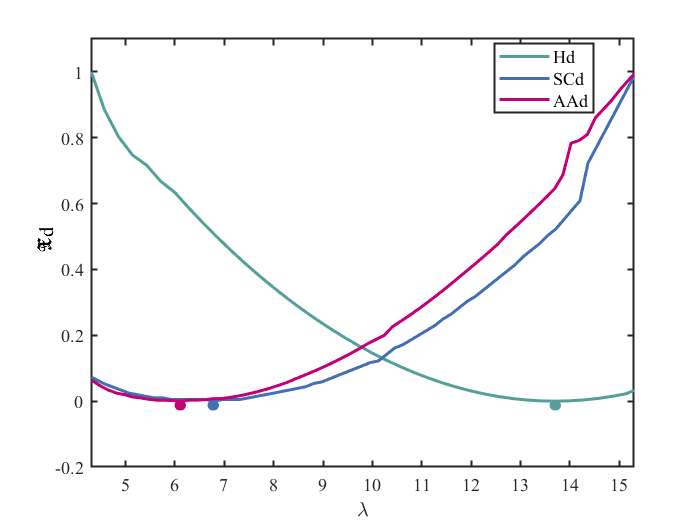}
\captionsetup{skip=1pt}
\caption{Different criteria measure for various $\lambda$ when $\hat{N}$ = 50 via PDE-based method by MFS in Example 3.}
\label{Figure 5}
\end{figure}
\vspace{-1cm} 
\begin{table}[htbp]
\vspace{3mm}
  \centering
  \captionsetup{skip=2pt}
  \caption{The volume of the Bunny model and related Hd by Kansa and MQ in Example 3}
    \begin{tabular}{c|c|c|c}
    \hline 
    \makebox[0.15\textwidth][c]{Method} &
    \makebox[0.15\textwidth][c]{$\hat{N}$} &
    \makebox[0.23\textwidth][c]{Volume} &
    \makebox[0.23\textwidth][c]{Hd} \\
    \hline \multirow{3}[0]{*}{Kansa} & 50    & 0.7533 &  0.0084 \\
                                  & 100   & 0.7542  & 0.0088 \\
                                  & 150   &0.7541  & 0.0088\\
    \hline \multirow{3}[0]{*}{MQ} & 50    &0.7544  & 0.0076 \\
                                  & 100   &0.7554 & 0.0083 \\
                                  & 150   & 0.7550 & 0.0090 \\
    \hline  
    \end{tabular}%
  \label{Table 33}%
\end{table}%

\newpage
\begin{figure}[htbp]
\centering
    \subfigure[(a) PDE-based method by MFS]{
	\includegraphics[width=10cm]{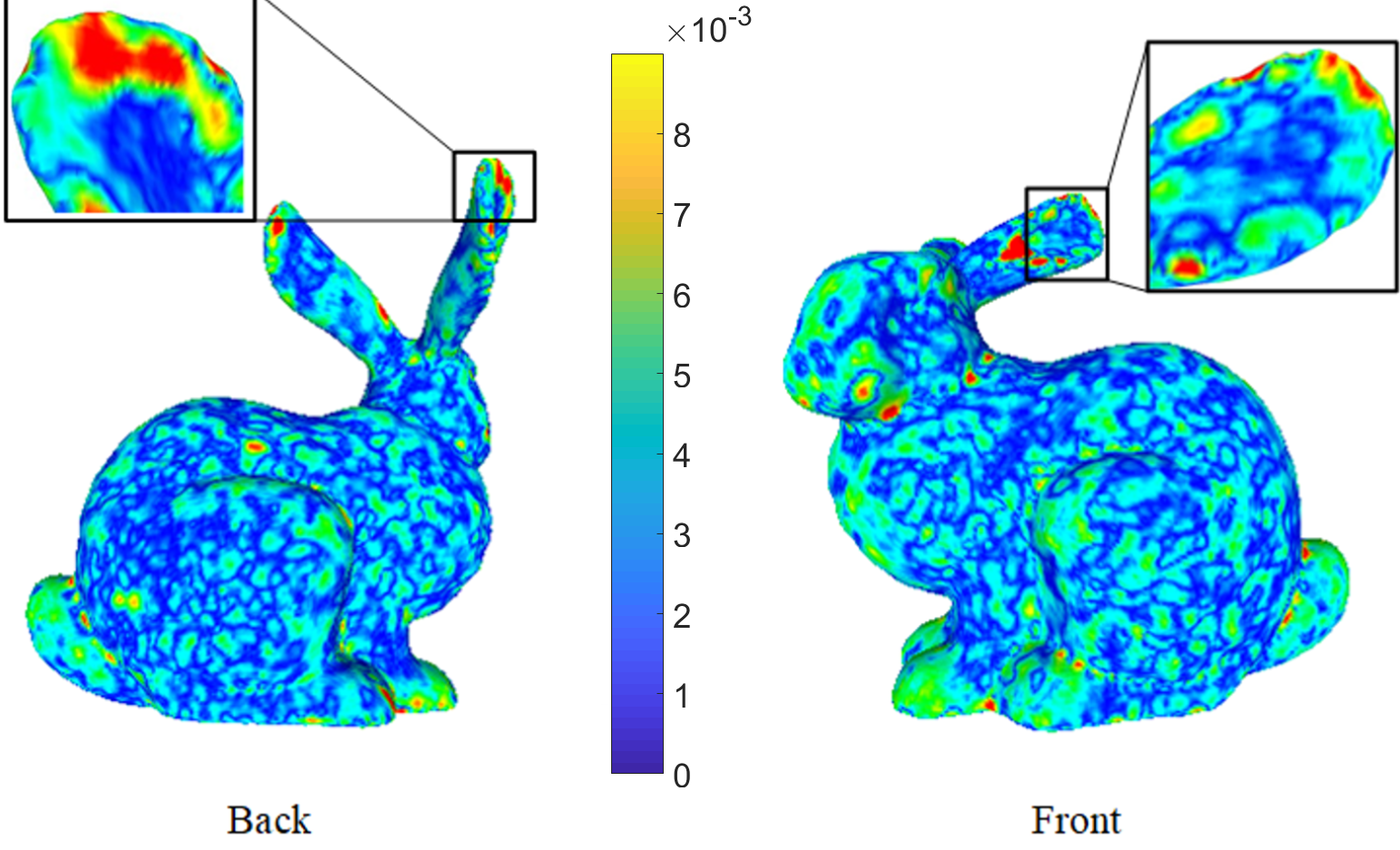} 
    \label{Fig.3(aa)}
    }	   
    \hspace{2mm}
	\subfigure[(b) PDE-based method by Kansa]{
	\includegraphics[width=10cm]{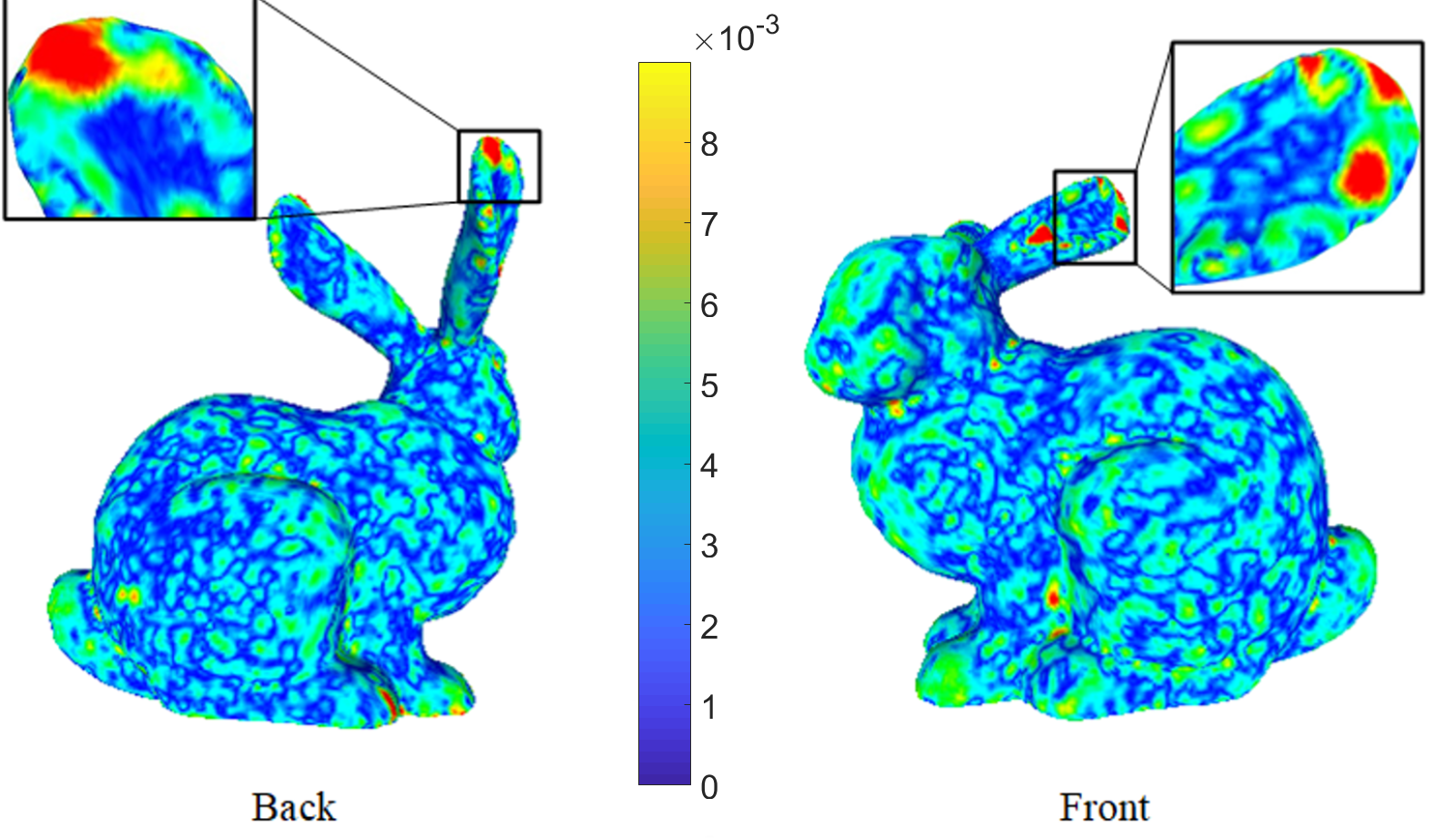} 
    \label{Fig.3(bbb)}
    }    
	\hspace{2mm}
	\subfigure[(c) Interpolation-based method by RBF-MQ]{
	\includegraphics[width=10cm]{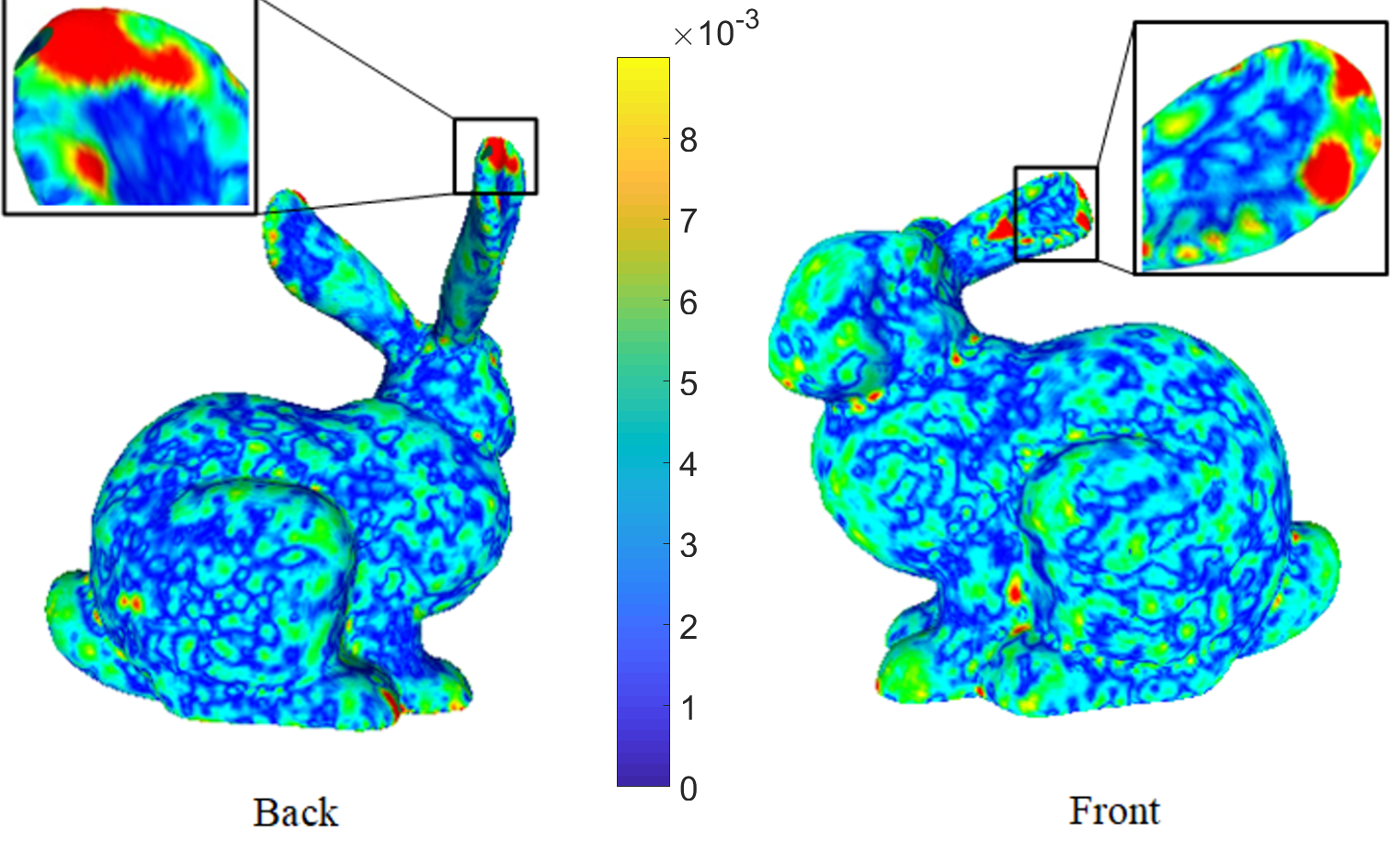} 
    \label{Fig.3(c)}
	}
 	\caption{The error distribution of surface reconstruction in Example 3.
  {
  \textcolor[rgb]{1.00,0.00,0.00}{
  (Plotted by MeshLab)
  }
  }
  }
\label{Figure Ex3_1}
\end {figure}

{
\newpage
\begin{figure}[htbp]
\centering 
\includegraphics[height=6.8cm,width=11cm]{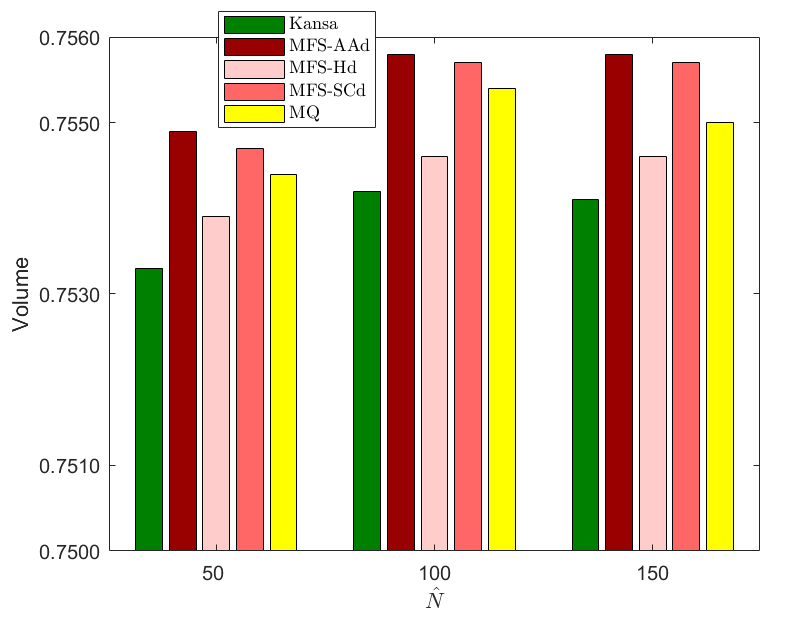}
\captionsetup{skip=1pt}
\caption{Volume data under different methods and criteria for Example 3.}
\label{Figure 5_1}
\end{figure}
\vspace{-0.7cm}
\subsection*{4.4 Example 4}
In this example, we consider the volume estimation of the dragon with 22998 boundary cloud points. 
\textcolor[rgb]{0.00,0.00,0.00}{
Figure \ref{Figure 9} indicates that the minimum Hd is 0.00583 when $\lambda^{*} = 28.2809 $ for MFS while the optimal ones found by SCd and AAd are  $\lambda^{*} = 26.4719 $ and  $\lambda^{*} = 24.6629 $.}
In Table \ref{Table 4_0}, the volume of the dragon with $\lambda^{*}$ found by Hd, SCd and AAd in MFS when $\hat{N}= 50,100,150$ are given respectively.
It can be seen that inconsistent optimal $\lambda^{*}$ in Figure \ref{Figure 9} 
will result in consistent volume estimates, with fluctuations within the range of $10^{-3}$.
 For the Kansa case, the minimum of Hd
 is 0.0052 when $\lambda = 95 $.
The volume of the dragon model and related Hd are shown in Table \ref{Table 4}.
The data in Figure \ref{Figure 6_1} reveals that while there are variations in volume estimates among the methods, the overall consistency in volume calculations is evident, with fluctuations generally within the range of $10^{-3}$.
This consistency across different algorithms and criteria underscores the robustness of the volume estimation process, ensuring reliable and accurate volume measurements regardless of the method employed.
From Figure \ref{Figure 6_1}, we can conclude that the volume is around 0.477 with oscillations $10^{-3}$, which is highly consistent with the results given by Meshlab (0.479).
The error distribution obtained by MFS, Kansa and MQ are shown in Figure \ref{Figure Ex4_1}, where we find MFS gives the best overall performance.
It shows that the main errors come from the dragon's head and hind claws.
The studied three methods almost have same level of accuracy 
in relatively flat areas and the main differences in the error distribution are emphasised in the red circle in the insets respectively.
\begin{table}[htbp]
  \centering
  \captionsetup{skip=1pt}
  \caption{Volumes with $\lambda^{*}$ given in Figure \ref{Figure 9} and related Hd, SCd and AAd in MFS for Example 4}
    \begin{tabular}{c|c|c|c}
    \hline 
    \makebox[0.1\textwidth][c]{$\hat{N}$} &
    \makebox[0.26\textwidth][c]{Volume (Hd)} & 
    \makebox[0.26\textwidth][c]{Volume (SCd)} &
    \makebox[0.26\textwidth][c]{Volume (AAd)} \\
    \hline 
    50   & 0.4771 ($4.3549e-03$) & 0.4772 ($2.1928e-06$) &0.4772 ($9.9574e-04$)\\ 
    100  & 0.4771 ($4.2836e-03$)& 0.4771 ($9.7902e-07	$)&0.4771 ($6.7294e-04$)\\  	
    150  & 0.4772 ($4.3811e-03$) & 0.4772 ($7.7332e-07$) & 0.4773 ($6.6572e-04$)\\ 	
    \hline  
    \end{tabular}%
  \label{Table 4_0}%
\end{table}%

\newpage
\begin{figure}[htbp]
\centering 
\includegraphics[height=7.8cm,width=9cm] {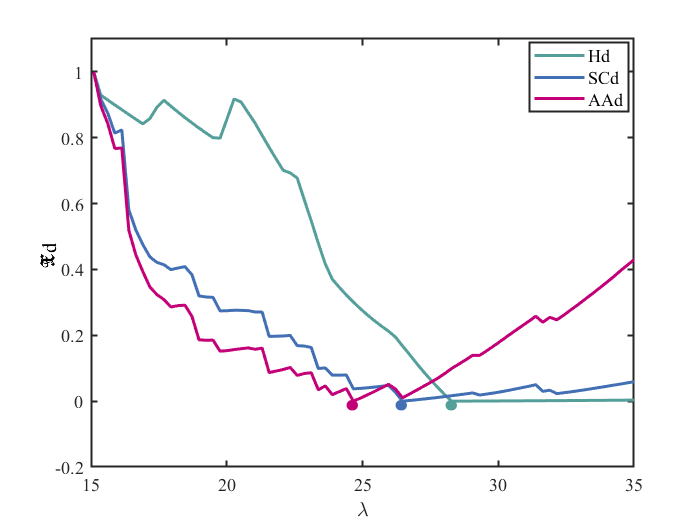}
\captionsetup{skip=1pt}
\caption{Different criteria for various $\lambda$ when $\hat{N}$ = 50 via MFS in Example 4.}
\label{Figure 9}
\end{figure}

\begin{table}[htbp]
\vspace{3mm}
  \centering
  \captionsetup{skip=1pt}
  \caption{The volume of the dragon model and related Hd by Kansa and MQ in Example 4}
    \begin{tabular}{c|c|c|c}
    \hline 
    \makebox[0.15\textwidth][c]{Method} &
    \makebox[0.15\textwidth][c]{$\hat{N}$} &
    \makebox[0.23\textwidth][c]{Volume} &
    \makebox[0.23\textwidth][c]{Hd} \\
    \hline \multirow{3}[0]{*}{Kansa} & 50  &0.4767  &0.0080   \\
                                   & 100   &0.4762  &0.0084   \\
                                   & 150   &0.4763   & 0.0084 \\
    \hline \multirow{3}[0]{*}{MQ} & 50    &0.4785  &0.0057   \\
                                   & 100   &0.4785   &0.0057  \\
                                   & 150   &0.4789  & 0.0051 \\
    \hline  
    \end{tabular}%
  \label{Table 4}%
\end{table}%

\vspace{0.5cm}
\begin{figure}[htbp]
\centering 
\includegraphics[height=6.8cm,width=8cm]{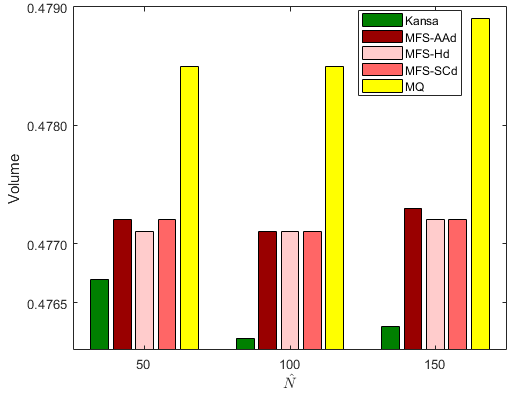}
\captionsetup{skip=1pt}
\caption{Volume data under different methods and criteria in Example 4.}
\label{Figure 6_1}
\end{figure}

\newgeometry{bottom=1.5cm}
{
\begin{figure}[htbp]
\centering
    \subfigure[(a) PDE-based method by MFS]{
	\includegraphics[width=12cm]{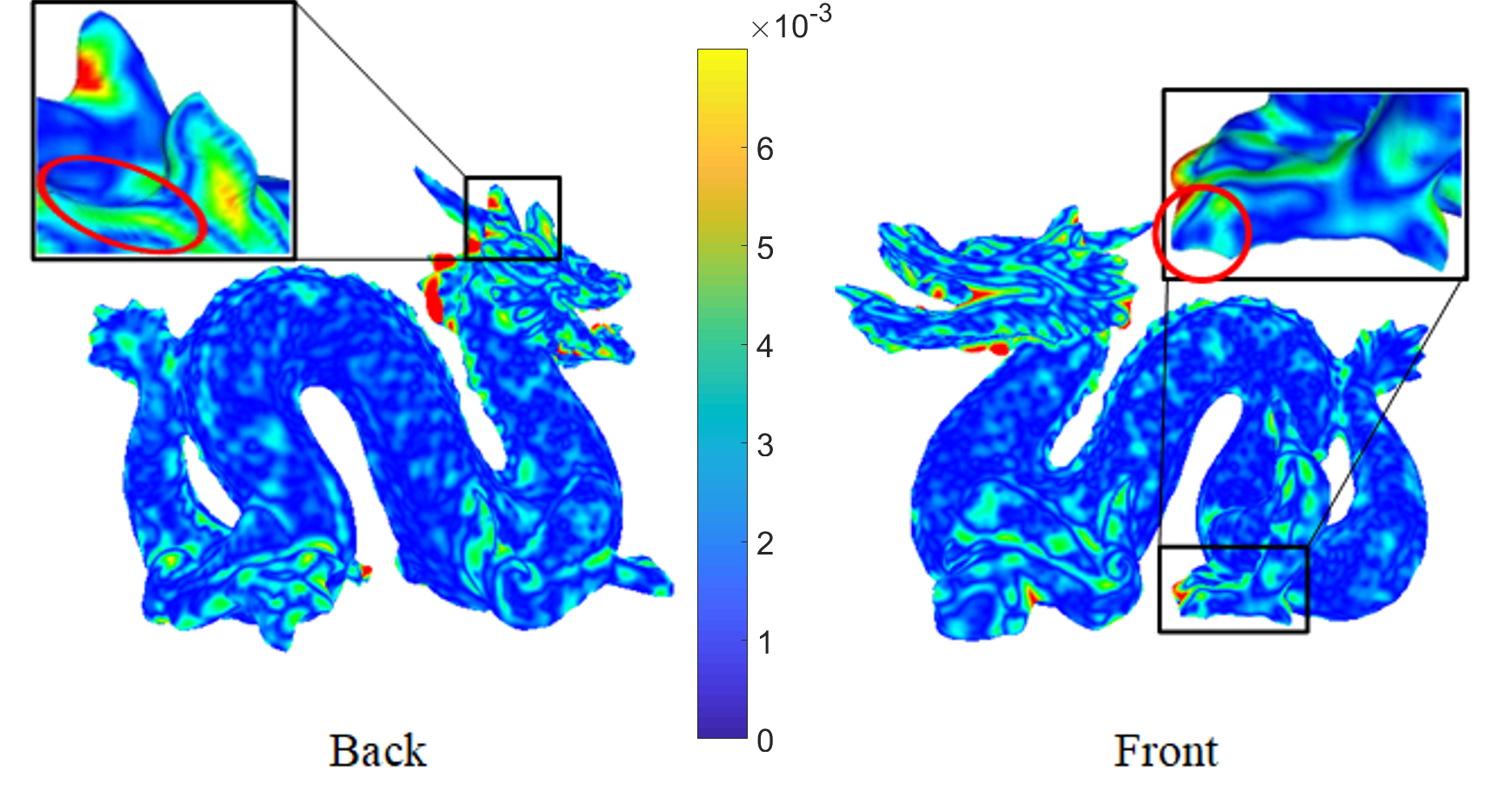} 
    \label{Fig.3(abb)}
    }	   
    \hspace{-5mm}
	\subfigure[(b) PDE-based method by Kansa]{
	\includegraphics[width=12cm]{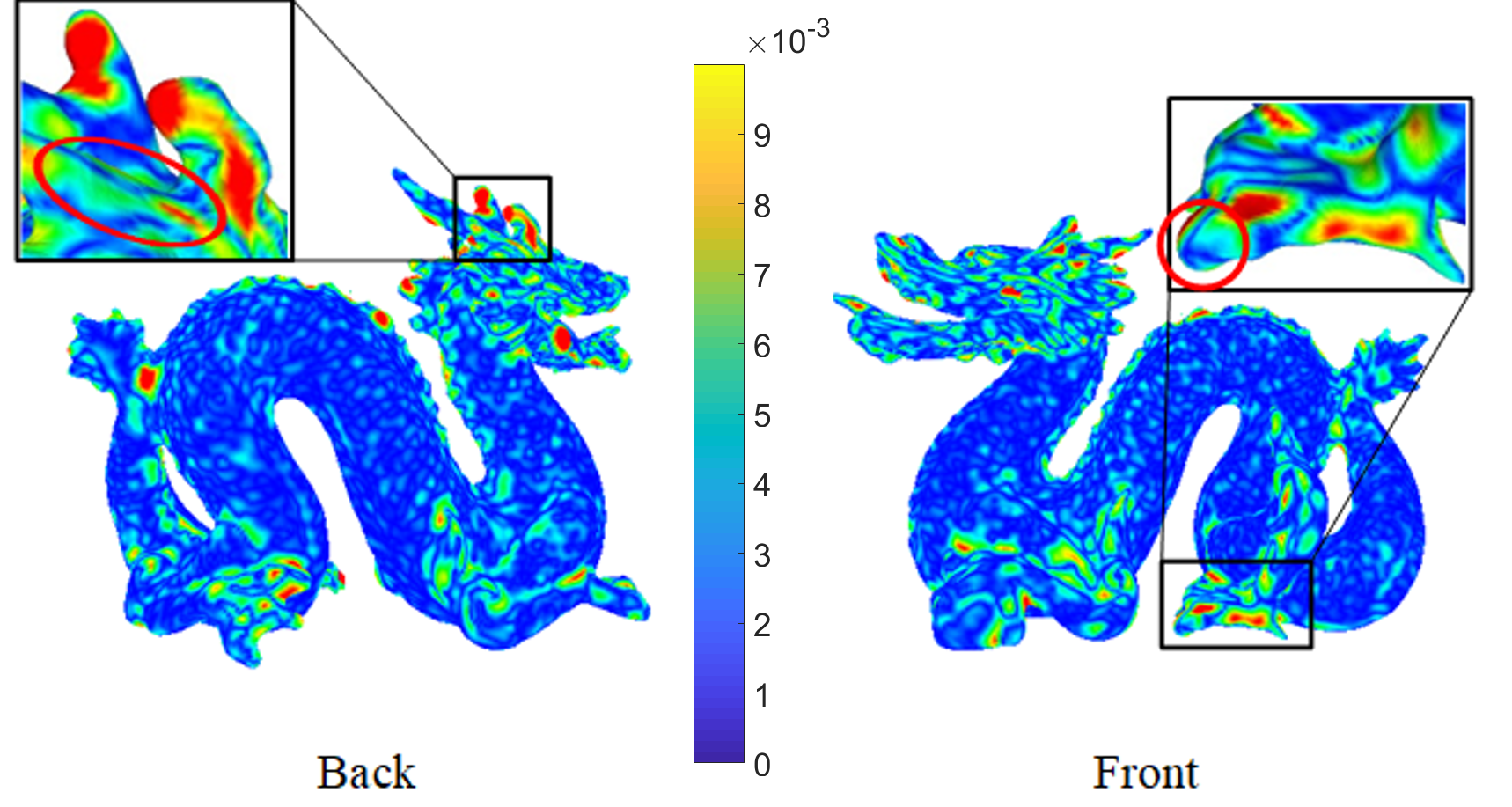} 
    \label{Fig.3(bb)}
    }    
	\hspace{-5mm}
	\subfigure[(c) Interpolation-based method by RBF-MQ]{
	\includegraphics[width=12cm]{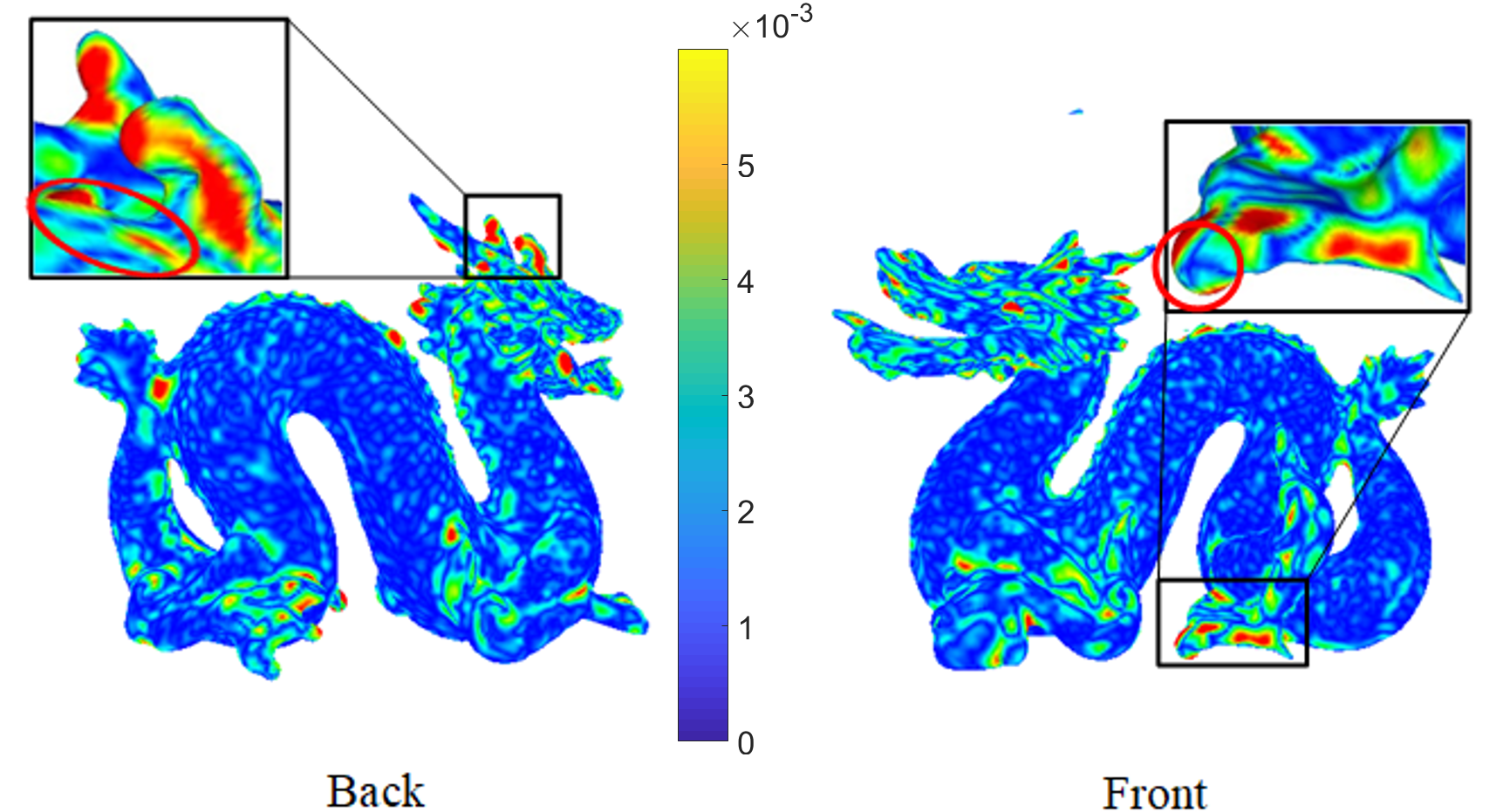} 
    \label{Fig.3(cb)}
	}
 	\caption{The error distribution of surface reconstruction in Example 4.
  {
  \textcolor[rgb]{1.00,0.00,0.00}{
  (Plotted by MeshLab)
  }
  }
  }
\label{Figure Ex4_1}
\end {figure}
}

\clearpage
\newpage
\section*{5. Conclusions}
In this paper, we study several approaches for the 3D implicit surface reconstruction from a set of scattered cloud data. 
Both the interpolation-based approach and PDE-based approach are considered.
By comparing above methods, it can be found that the PDE-based method using MFS will result in a smaller system matrix as well as a higher efficiency.
The overall performance of PDE-based approach with MFS is the best.
Besides, the Hausdorff distance, Symmetric Chamfer distance, Absolute Average distance are employed to determine the free parameter $\lambda$ in proposed PDE models in MFS and Kansa.
Based on the optimal reconstructed models, the volume estimation strategy is proposed further. 
This strategy effectively estimates the volume of 3D models, even those with small holes. Compared with traditional tools like MeshLab, which cannot calculate the volume of non-'watertight' models, our method is particularly suitable for specialized 3D models, such as incomplete vascular plaques often encountered in medical scenarios.
Plenty of numerical examples show that the varying optimal $\lambda^{*}$ values determined by the Hd, SCd, and AAd criteria yield consistent volume estimation. This uniformity highlights the resilience of the surface reconstruction and volume estimation process, guaranteeing dependable and precise outcomes regardless of the criterion applied.

\par
\section*{Acknowledgements}
\indent The author is grateful for the support of the National Youth Science Foundation of China (No.12102283) and Youth Science and Technology Research Foundation of Shanxi Province (No. 202103021223059 and  No.20210302124159).

\section*{References}
\bibliography{References}
\end{document}